\documentclass[12pt, journal, draftclsnofoot, onecolumn]{IEEEtran}
\ifCLASSOPTIONcompsoc
\else
\fi

\ifCLASSINFOpdf
\else
\fi
\usepackage[parfill]{parskip}    
\usepackage{amsmath, amssymb, latexsym, enumerate, mathrsfs}
\usepackage{graphicx}
\usepackage{hyperref}
\usepackage{graphicx}
\usepackage{array}
\usepackage{times}

\usepackage{moreverb}
\usepackage{marvosym}
\usepackage{color,soul}
\definecolor{lightblue}{rgb}{.90,.95,1}
\sethlcolor{yellow}

\usepackage{hyperref}

\newtheorem{lemma}{\bf Lemma}
\newtheorem{proposition}{\bf Proposition}

\newtheorem{remark}{Remark}

\newtheorem{assumption}{\bf Assumption}

\def\qed{\hfill $\Box$}

\begin{document}
%
\title{Large Deviations for the Reliability Assessment of Redundant Multi-Channel Systems}

\author{Getachew K. Befekadu and Panos J. Antsaklis
\thanks{}
\IEEEcompsocitemizethanks{\IEEEcompsocthanksitem G. K. Befekadu is with the Department
of Electrical Engineering, University of Notre Dame, Notre Dame, IN 46556, USA.\protect\\
E-mail: gbefekadu1@nd.edu
\IEEEcompsocthanksitem P. J. Antsaklis is with the Department
of Electrical Engineering, University of Notre Dame, Notre Dame, IN 46556, USA.\protect\\
E-mail: antsaklis.1@nd.edu}}

\markboth{}%
{Shell \MakeLowercase{\textit{et al.}}: Bare Advanced Demo of IEEEtran.cls for Journals}
\IEEEcompsoctitleabstractindextext{%
\begin{abstract}
In this paper, we are concerned with the reliability assessment of redundant multi-channel systems having multiple controllers with overlapping functionality -- where all controllers are required to respond optimally to the non-faulty controllers so as to ensure or maintain some system properties. In particular, for such redundant systems with small random perturbation, we study the relationships between the exit probabilities with which the state-trajectories exit from a given bounded open domain and the value functions corresponding to a family of stochastic exit-time control problems on the boundary of the given domain. Moreover, as the random perturbation vanishes, such relationships provide useful information concerning the reliability of the redundant multi-channel systems arising from the large deviations problem in connection with the asymptotic estimates of exit probabilities with respect to some portions of the boundary of the given domain. Finally, we briefly comment on the implication of our results on a co-design technique using multi-objective optimization frameworks for evaluating the performance of the redundant multi-channel systems.
\end{abstract}

}

\maketitle

\IEEEdisplaynotcompsoctitleabstractindextext

%
\IEEEpeerreviewmaketitle

\section{Introduction} \label{S1}
In this paper, we are concerned with the reliability assessment of redundant multi-channel systems having multiple controllers with overlapping functionality. Specifically, we consider a redundant system with multi-controller configurations -- where all controllers are required to respond optimally, in the sense of best-response correspondence, i.e., a {\it reliable-by-design} requirement, to the non-faulty controllers so as to ensure or maintain some system properties. Here we are mainly interested in a systematic understanding of the relationships between the exit probabilities with which the state-trajectories exit from a given bounded open domain and the value functions corresponding to a family of stochastic exit-time control problems on the boundary of the given domain. As a consequence of such relationships, we obtain useful information concerning the reliability of the redundant multi-channel systems arising from the large deviations problem as the random perturbation vanishes in connection with the asymptotic estimates of exit probabilities with respect to some portions of the boundary of the given domain. Moreover, we also comment on the implication of our results on a co-design technique using multi-objective optimization frameworks for either evaluating the performance or finding an appropriate set of redundant controllers for the multi-channel systems with respect to some prescribed portions of the boundary of the given domain.

It is worth mentioning that some interesting studies on the exit probabilities for the dynamical systems with small random perturbation have been reported in literature (see, e.g., \cite{VenFre70}, \cite{FreWe84}, \cite{Kif90} and \cite{DupKu86} in the context of large deviations; see \cite{Day86}, \cite{EvaIsh85}, \cite{Fle78} and \cite{FleTs81} in connection with optimal stochastic control problems; and see \cite{Day86} or \cite{MatSc77} via asymptotic expansions approach). Note that the rationale behind our framework follows, in some sense, the settings of these papers -- where we establish a connection between the asymptotic estimates of the exit probabilities and the stochastic exit-time control problems on some portions of the boundary of the given domain. However, to our knowledge, such a connection has not been addressed in the context of multi-channel systems with multi-controller configurations having ``overlapping or backing-up" functionality, and it is important because it provides a framework that shows how the asymptotic estimates on the exit probabilities can be systematically used to obtain useful information concerning the reliability of the redundant multi-channel systems.\footnote{In this paper, our intent is to provide a theoretical result, rather than considering a specific numerical problem or application.}

The rest of the paper is organized as follows. In Section~\ref{S2}, we present some preliminary results that are useful for our main results. In Section~\ref{S3}, we briefly discuss a family of boundary value problems for the multi-channel system in the presence of small random perturbations. In Section~\ref{S4}, we provide asymptotic estimates on the exit probabilities on the positions of the state-trajectories at the first time of their exit through some portions of the boundary of the given domain. This section also provides connections between such asymptotic estimates and the value functions corresponding to a family of stochastic exit-time control problems on the boundary of the given domain. Moreover, we use such asymptotic estimates to obtain useful information on the reliability of the redundant multi-channel systems. Finally, Section~\ref{S5} provides some further remarks.

\section{Preliminaries} \label{S2}
Consider the following continuous-time multi-channel system 
\begin{align} 
 \dot{x}(t) &= A x(t) +  \sum\nolimits_{i=1}^n B_i u_i(t), \quad x(0)=x_0,  \label{Eq1}
\end{align}
where $A \in \mathbb{R}^{d \times d}$, $B_i \in \mathbb{R}^{d \times r_i}$, $x(t) \in \mathbb{R}^{d}$ is the state of the system, $u_i(t) \in \mathbb{R}^{r_i}$ is the control input to the $i$th\,-\,channel in the system. 

In what follows, we consider a particular class of stabilizing state-feedbacks that satisfies\footnote{$\operatorname{Sp}(A)$ denotes the spectrum of a matrix $A \in \mathbb{R}^{d \times d}$, i.e., $\operatorname{Sp}(A)=\bigl\{ s \in \mathbb{C}\, \vert \, \operatorname{rank}(A-sI) < d\bigr\}$.}
\begin{align}
  \mathcal{K} \subseteq \biggl\{\underbrace{\bigl(K_1, K_2, \ldots, K_n\bigr)}_{\substack{\triangleq K}} \in \prod\nolimits_{i=1}^n \mathbb{R}^{r_i \times d} \biggm\lvert  \operatorname{Sp}\Bigl(A + \sum\nolimits_{i=1}^n B_{i} K_i\Bigr) \subset \mathbb{C}^{-} \quad \quad\quad\quad & \notag \\
    \operatorname{Sp}\Bigl(A + \sum\nolimits_{i \neq j} B_{i} K_i\Bigr) \subset \mathbb{C}^{-}, \,\,  j = 1,2, \ldots, n \Biggr\}, & \label{Eq2}
\end{align}

\begin{remark} \label{R1}
We remark that the above class of state-feedbacks is useful for maintaining the stability of the closed-loop system both when all of the controllers work together, i.e., $\bigl(A + \sum\nolimits_{i=1}^n B_{i} K_i\bigr)$, as well as when there is a single-channel controller failure in the system, i.e., $\bigl(A + \sum\nolimits_{i \neq j} B_{i} K_i\bigr)$ for $j \in \{1,2, \ldots, n\}$. Moreover, such a class of state-feedbacks falls within the redundant/passive fault tolerant multi-controller configurations with overlapping functionality (see, e.g., \cite{BefGA14} or \cite{FujBe09} for such a {\it reliable-by-design} requirement in multi-channel systems).   
\end{remark}

Consider the following family of stochastic differential equations
\begin{align} 
 dx^{\epsilon, 0}(t) = \Bigl(A + \sum\nolimits_{i=1}^ n B_i K_i \Bigr)x^{\epsilon, 0}(t) dt + \sqrt{\epsilon} \,\sigma\bigl(x^{\epsilon, 0}(t)\bigr) dW(t), \,\, x^{\epsilon, 0}(0)=x_0 \label{Eq3} 
\end{align}
and
\begin{align} 
 dx^{\epsilon, j}(t) = \Bigl(A + \sum\nolimits_{i\neq j} B_i K_i \Bigr)x^{\epsilon, j}(t) dt + \sqrt{\epsilon} \,\sigma\bigl(x^{\epsilon, j}(t)\bigr) dW(t),\,\, x^{\epsilon, j}(0)=x_0, \,\, j = 1, 2, \ldots, n,  \label{Eq4} 
 \end{align}
where $\sigma \in \mathbb{R}^{d \times d}$ is a diffusion term, $W$ (with $W(0)=0$) is a $d$-dimensional Wiener process and $\epsilon$ is a small positive number, which represents the level of perturbation in the system. 

Let $\Omega \subset \mathbb{R}^{d}$ be a bounded open domain with smooth boundary (i.e., $\partial \Omega$ is a manifold of class $C^2$). Then, the second-order elliptic differential operators $\mathcal{L}^{(\epsilon, j)}$ that correspond to the above family of stochastic differential systems are given by
\begin{align}
 \mathcal{L}^{(\epsilon, 0)}f^{(j)}(x)= \frac{\epsilon}{2} \sum_{i,k=1}^d a_{ik}(x) \frac{\partial^2f^{(j)}(x)}{\partial x_i \partial x_k} + \Bigl \langle \Bigl(A + \sum\nolimits_{i=1}^n B_i K_i \Bigr)x,\, \triangledown f^{(j)}(x) \Bigr \rangle \label{Eq5}
\end{align}
and
\begin{align}
 \mathcal{L}^{(\epsilon, j)}f^{(j)}(x)= \frac{\epsilon}{2} \sum_{i,k=1}^d a_{ik}(x) \frac{\partial^2f^{(j)}(x)}{\partial x_i \partial x_k} + \Bigl \langle \Bigl(A + \sum\nolimits_{i \neq j} B_i K_i \Bigr)x,\, \triangledown f^{(j)}(x) \Bigr \rangle, \quad j = 1, 2, \ldots, n, \label{Eq6}
\end{align}
where $f^{(j)}(x) \in C^2(\Omega) \cap C^1(\bar{\Omega})$ and $\triangledown f^{(j)}(x)$ is the gradient of $f^{(j)}(x)$. Further, we assume that the matrix $a(x)=\sigma(x)\sigma^T(x)$ is nonnegative definite and $\sigma(x)$ satisfies a global Lipschitz condition. Hence, the operators $\mathcal{L}^{(\epsilon, j)}$ are uniformly elliptic for fixed $\epsilon > 0$. 

Let $C\bigl([0,\infty),\mathbb{R}^{d}\bigr)$ denote the space of continuous functions from $[0,\infty)$ to $\mathbb{R}^{d}$, and let $\mathcal{H}^{1}[0, T]$ be the space of all $\varphi \in C\bigl([0,\infty),\mathbb{R}^{d}\bigr)$ such that $\varphi(t)$ is absolutely continuous and $\int_0^T\vert \dot{\varphi}(t)\vert dt < \infty$ for each $T > 0$. Let us associate portions of the boundary $\Gamma_j \subset \partial \Omega$ for $j=0,1, \ldots, n$, with different operating conditions of the redundant system (for example, $\Gamma_0$ with the nominal operating condition and $\Gamma_j$ for $j = 1,2, \ldots, n$, with any of single-channel failures in the system). Then, as $\epsilon \rightarrow 0$, we investigate the behavior of the solutions for the second-order elliptic equations corresponding to a family of boundary value problems (see equation~\eqref{Eq12}) with respect to those portions of the boundary $\Gamma_j \subset \partial \Omega$ for $j=0,1, \ldots, n$ (see \cite{Lev50} for discussions on the first boundary value problem with small parameter). In general, such asymptotic estimates involve finding a family of minimum functionals $I_j\bigl(\varphi, \tau_{j}\bigr)$, i.e., 
\begin{align}
 I_j\bigl(\varphi, \tau_{j}\bigr) = \inf_{\varphi \in \mathcal{H}^{1}[0, T], \tau_{j} \ge 0} \int_{0}^{\tau_{j} \wedge T} L_j\big(\varphi(t),\dot{\varphi}(t)\big)dt,\quad  j=0, 1, \ldots, n, \label{Eq7}
\end{align}
where the infimum is taken among all $\varphi \in \mathcal{H}^{1}[0, T]$ and $\tau_{j} \ge 0$ (where $\tau_{j}$ is the first exit-time of $x^{\epsilon, j}(t)$ from $\Omega$\footnote{$\tau_{j} = \inf \bigl\{ t \,\vert \,x^{\epsilon, j}(t) \notin \Omega \bigr\}$.}) such that $\varphi(0)=x_0$ and $\varphi(t) \in \bar{\Omega}$ for $t \in [0, \tau_{j} \wedge T]$, with
\begin{align*}
 L_0 \bigl(\varphi(t), \dot{\varphi}(t)\bigr) = \frac{1}{2} \biggm\Vert\dot{\varphi}(t) - \Bigl(A + \sum\nolimits_{i=1}^n B_i K_i \Bigr)\varphi(t)\biggm\Vert_{\bigm[a(\varphi(t))\bigm]^{-1}}^2 
\end{align*}
and
\begin{align*}
 L_j \bigl(\varphi(t), \dot{\varphi}(t)\bigr) = \frac{1}{2} \biggm\Vert\dot{\varphi}(t) - \Bigl(A + \sum\nolimits_{i \neq j} B_i K_i \Bigr)\varphi(t)\biggm\Vert_{\bigm[a(\varphi(t))\bigm]^{-1}}^2, \quad j=1,2, \ldots, n,
\end{align*}
where $a(\varphi(t)) = \sigma(\varphi(t)) \sigma^T(\varphi(t))$.\footnote{$\Vert x\Vert_{P}^2 \triangleq x^T P x, \quad x \in \Omega$.} 

Furthermore, if we let $\varphi(0)=x_0$ in the domain $\Omega$, then, for any $\Gamma_j \subset \partial \Omega$, the infimum in \eqref{Eq7}, when subjected to an additional condition $\varphi(\tau_{j}) \in \Gamma_j$, will attain 
\begin{align}
 I_j\bigl(x_0, \Gamma_j \bigr) = -\lim_{\epsilon \rightarrow 0} \epsilon\,\log \mathbb{P}_{x_0}^{(\epsilon, j)}\Bigl(x^{\epsilon, j}(\tau_{j}) \in \Gamma_j \Bigr), \quad j=0,1, \ldots, n,
\label{Eq8}
\end{align}
which implicitly depends on the initial condition $x_0$ and the boundary $\Gamma_j$. Note that such information, which is based on \eqref{Eq8}, will be used to identify the exit positions on the boundary $\Gamma_j \subset \partial \Omega$ under additional assumptions on the behavior of the state-trajectories of the unperturbed systems $x^{\epsilon, j}$, when $\epsilon=0$, as $t \rightarrow \infty$.

In Sections~\ref{S3} and \ref{S4}, we establish relationships between the exit probabilities with which the state-trajectories exit from a given bounded open domain and the value functions corresponding to a family of stochastic exit-time control problems on the boundary of the given domain. More specifically, we provide asymptotic estimates for the exit probabilities on the positions of the state-trajectories $x^{\epsilon, j}(t)$, for each $j=0,1, \ldots, n$, at the first time of their exit from a bounded open domain $\Omega \subset \mathbb{R}^{d}$ (i.e., estimating bounds on the exit probabilities of the state-trajectories $x^{\epsilon, j}(t)$ from the given domain $\Omega$ through a portion or section of the given boundary $\Gamma_j \subset \partial \Omega$ (see Proposition~\ref{P1})). Such asymptotic estimates (i.e., the asymptotic estimates on the exit probabilities $\mathbb{P}_{x_0}^{(\epsilon, j)}\bigl(x^{\epsilon, j}(\tau_{j}) \in \Gamma_j\bigr)$, as $\epsilon \rightarrow 0$, conditioned on the initial point $x_0 \in \Omega$) can be linked to finding probabilities for the state-trajectories $x^{\epsilon, j}(t)$ that do not deviate by more than $\delta$ from a smooth function $\varphi \in \mathcal{H}^{1}[0, T]$ during the time $t \in [0,\tau_{j} \wedge T]$.\footnote{Note that the behavior of $-\epsilon \log \mathbb{P}_{x_0}^{(\epsilon, j)}\bigl(x^{\epsilon, j}(\tau_{j}) \in \Gamma_j\bigr)$, as $\epsilon \rightarrow 0$, is defined by the large deviations of the state-trajectories from their typical behavior.} Moreover, for small $\delta > 0$, the exit probabilities $\mathbb{P}_{x_0}^{(\epsilon, j)}\bigl(x^{\epsilon, j}(\tau_{j}) \in \Gamma_j\bigr)$ will have forms $\exp\bigl(-\dfrac{1}{\epsilon}I_j(x_0, \Gamma_j)\bigr)$, where $I_j(x_0, \Gamma_j)$ is a non-negative functional of $\varphi \in \mathcal{H}^{1}[0, T]$ (see Proposition~\ref{P2}). 

\section{Boundary value problem} \label{S3}
For the family of stochastic differential equations in \eqref{Eq3} and \eqref{Eq4}, consider the following family of boundary value problems
\begin{align}
\left.\begin{array}{l}
 \mathcal{L}^{(\epsilon, j)} f^{(j)}(x) = 0 \quad \text{in} \quad \Omega\\
 \quad\quad f^{(j)}(x) = \mathbb{E}_{x_0}^{(\epsilon, j)} \Bigl(\exp\Bigl(-\frac{1} {\epsilon} \Phi_j\bigl(x\bigr) \Bigr) \Bigr) \quad \text{on} \quad \partial\Omega\\
 \quad\quad\quad\quad  j = 0, 1, \ldots, n
\end{array} \right\} \label{Eq9}
\end{align}
where $\Phi_j$ is class $C^2$ function, with $\Phi_j \ge 0$. Then, there exists a set of unique solutions $f^{(j)}(x) \in C^2(\Omega) \cap C^1(\bar{\Omega})$ such that 
\begin{align}
  f^{(j)}(x) = \mathbb{E}_{x_0}^{(\epsilon, j)} \Bigl(\exp\Bigl(-\frac{1}{\epsilon} \Phi_j \bigl(x^{\epsilon, j}(\tau_{j}) \bigr) \Bigr) \Bigr), \label{Eq10}
 \end{align}
where $\tau_{j}$ is the exit-time of $x^{\epsilon, j}(t)$ from the domain $\Omega$. Note that if we further introduce the following logarithmic transformation (see, e.g., \cite{EvaIsh85} or \cite{Fle78} for such logarithmic connections between large deviations and stochastic optimization problems)
\begin{align}
  J_{\Phi}^{(\epsilon, j)}(x) &= - \epsilon \log \Bigl(f^{(j)}(x)\Bigr), \notag \\
                            &= - \epsilon \log \Bigl(\mathbb{E}_{x_0}^{(\epsilon, j)} \Bigl(\exp\Bigl(-\frac{1}{\epsilon} \Phi_j \bigl(x^{\epsilon, j}(\tau_{j}) \bigr) \Bigr) \Bigr)\Bigr), \quad j=0,1, \ldots, n. \label{Eq11}
\end{align}
Then, $J_{\Phi}^{(\epsilon, j)}(x)$ satisfies the following second-order elliptic differential equation
\begin{align}
0 = \frac{\epsilon}{2} \sum_{i,k=1}^d a_{ik}(x) \frac{\partial^2 J_{\Phi}^{(\epsilon, j)}(x)}{\partial x_i \partial x_k} + H_j \bigl(x, \triangledown J_{\Phi}^{(\epsilon, j)}(x)\bigr) \quad \text{in} \quad \Omega, & \notag \\
 j=0,1, \ldots, n, & \label{Eq12}
\end{align}
where
\begin{align*}
H_0 \bigl(x, \triangledown J_{\Phi}^{(\epsilon, 0)}(x)\bigr) = \frac{1}{2} \biggl\Vert \triangledown J_{\Phi}^{(\epsilon, 0)}(x)\biggr\Vert_{\bigm[a(\varphi(t))\bigm]^{-1}}^2 + \biggl\langle \triangledown J_{\Phi}^{(\epsilon, 0)}(x),\, \bigl(A + \sum\nolimits_{i=1}^n B_i K_i \bigr)x \biggr\rangle
\end{align*}
and 
\begin{align*}
H_j \bigl(x, \triangledown J_{\Phi}^{(\epsilon, j)}(x)\bigr) = \frac{1}{2} \biggl\Vert \triangledown J_{\Phi}^{(\epsilon, j)}(x)\biggr\Vert_{\bigm[a(\varphi(t))\bigm]^{-1}}^2 + \biggl\langle \triangledown J_{\Phi}^{(\epsilon, j)}(x),\, \bigl(A + \sum\nolimits_{i \neq j} B_i K_i \bigr)x \biggr\rangle, \\
 \,\, j=1,2,\ldots, n.
\end{align*}
Further, note that there is a duality between $H_j \bigl(x,\, \cdot\bigr)$ and $L_j \bigl(x,\, \cdot\bigr)$, for each $j=0, 1, \ldots, n$, such that
\begin{align}
H_j \bigl(x, \triangledown J_{\Phi}^{(\epsilon, j)}(x)\bigr) = \inf_{\upsilon} \biggl \{L_j \bigl(x, \triangledown J_{\Phi}^{(\epsilon, j)}(x)\bigr) + \Bigl\langle \triangledown J_{\Phi}^{(\epsilon, j)}(x),\, \upsilon \Bigr\rangle \biggr\}. \label{Eq13}
\end{align}
Hence, it is easy to see that $J_{\Phi}^{(\epsilon, j)}(x)$ is a solution in class $C^2(\Omega) \cap C^1(\bar{\Omega})$, with $J_{\Phi}^{(\epsilon, j)}=\Phi_j$ on $\partial \Omega$, to the dynamic programming in \eqref{Eq12}, where the latter is associated with the following stochastic exit-time control problem 
\begin{align}
 J_{\Phi}^{(\epsilon, j)}(x_0, \upsilon^{(j)}) =  \mathbb{E}_{x_0}^{(\epsilon, j)}\left\{ \int_{0}^{\tau_{j} \wedge T} L_j\bigl(\eta^{(j)}(t),\upsilon^{(j)}(t)\bigr)dt + \Phi_j\bigl(\eta^{(j)}(\tau_{j})\bigr) \right\},\quad j=0, 1, \ldots, n, \label{Eq14}
\end{align}
and $\eta^{(j)}(t)$ satisfies the following stochastic differential equation
\begin{align}
 d \eta^{(j)}(t) = \upsilon^{(j)}(t) dt + \sqrt{\epsilon} \,\sigma\bigl(\eta^{(j)}(t)\bigr) dW(t) \label{Eq15}
\end{align}
for $j=0, 1, \ldots, n$ (see, e.g., \cite{DavVar73}).

In the following section, i.e., Section~\ref{S4}, we exploit this formalism to prove the asymptotic bounds (cf. Proposition~\ref{P1}) for the exit probabilities on the position of state-trajectories at the first time of their exit from the given portion or section of the boundary of the domain $\Omega$ (cf. Proposition~\ref{P2}).  
  
\section{Main results} \label{S4}
In this section, we present our main results, i.e., the asymptotic estimates bounds on the exit probabilities of the state-trajectories $x^{\epsilon, j}(t)$ from the given domain $\Omega$ through the given portion (or section) of the boundary $\Gamma_j \subset \partial \Omega$ for $j=0, 1, \ldots, n$. Further, for $\Gamma_j \subset \partial \Omega$ for $j=0, 1, \ldots, n$, and $x_0 \in \Omega$, let 
\begin{align}
 q^{(\epsilon,j)}\bigl(x_0, \Gamma_j\bigr) &= \mathbb{P}_{x_0}^{(\epsilon, j)}\bigm(x^{\epsilon, j}(\tau_{j}) \in \Gamma_j \bigm), \label{Eq16}
\end{align}
and
\begin{align}
 I_j\bigl(x_0, \Gamma_j\bigr) &= -\lim_{\epsilon \rightarrow 0} \epsilon\,\log \mathbb{P}_{x_0}^{(\epsilon, j)}\bigm(x^{\epsilon, j}(\tau_{j}) \in \Gamma_j \bigm), \label{Eq17}
\end{align}
where $\tau_{j}$ is the first exit-time of $x^{\epsilon, j}(t)$ from the domain $\Omega$. Moreover, let
\begin{align}
 I_j\bigl(\varphi, \tau_{j}\bigr) = \inf_{\varphi \in \mathcal{H}^{1}[0, T], \tau_{j} \ge 0} \int_{0}^{\tau_{j} \wedge T} L_j\big(\varphi(t),\dot{\varphi}(t)\big)dt,\quad j=0, 1, \ldots, n, \label{Eq18}
\end{align}
where the infimum is taken among all $\varphi \in \mathcal{H}^{1}[0, T]$ and $\tau_{j} \ge 0$ such that $\varphi(0)=x_0$, $\varphi(t) \in \bar{\Omega}$ for $t \in [0, \tau_{j} \wedge T]$ and $\varphi(\tau_{j}) \in \Gamma_j$. Then, we have
\begin{align}
 I_j\bigl(x_0, \Gamma_j\bigr) = I_j\bigl(x_0, \bar{\Gamma}_j\bigr), \quad j=0, 1, \ldots, n. \label{Eq19}
\end{align}

\begin{remark} \label{R3}
Note that the functional 
\begin{align*}
\int_{0}^{\tau_{j} \wedge T} L_j\big(\varphi(t),\dot{\varphi}(t)\big)dt, \quad j=0, 1, \ldots, n,
\end{align*}
is lower semicontinuous with respect to $\varphi$ and $\tau_{j} \wedge T$. Furthermore, the set level $\Psi=\bigl\{\varphi \in \mathcal{H}^{1}[0, T] \,\bigl\vert \,I_j(\varphi, \tau_{j}) \le \alpha\bigr\}$ is a compact subset of $\mathcal{H}^{1}[0, T]$ for every $\alpha \ge 0$ and $\varphi(0)=x_0 \in \Omega$. Hence, the infimum in \eqref{Eq18} attains a minimum on $\Psi$ (see, e.g., \cite[pp.\,332,\, Corollary~1.4]{Fre06} or \cite{Fle78}).
\end{remark}

Next, we introduce the following assumption about the domain $\Omega$, which is useful in the sequel.
\begin{assumption} \label{A1} 
If $\varphi \in \mathcal{H}^{1}[0, T]$ and $\varphi(t) \in \bar{\Omega}$ for all $t \ge 0$, then $\int_{0}^{T} L_j\big(\varphi(t),\dot{\varphi}(t)\big)dt \rightarrow  +\infty$ as $T \rightarrow \infty$ for each $j=0, 1, \ldots, n$.
\end{assumption}

Consider again the stochastic control problem in \eqref{Eq14} (together with equation~\eqref{Eq15}). Suppose that  $\Phi_M$ (with $\Phi_M \ge 0$) is class $C^2$ such that $\Phi_M(x) \rightarrow +\infty$ as $M \rightarrow \infty$ uniformly on any compact subset of $\Omega \setminus \bar{\Gamma}$ and $\Phi_M(x)$ on $\bar{\Gamma_j}$ for $j=0,1, \ldots, n$. Further, if we let $J_{\Phi}^{(\epsilon, j)}(x) = J_{\Phi_M}^{(\epsilon, j)}(x)$, when  $\Phi_j = \Phi_M$, then we have the following lemma.

\begin{lemma} \label{L1}
Suppose that Assumption~\ref{A1} holds, then we have
\begin{align}
 \liminf_{\substack{M \rightarrow \infty \\ x \rightarrow x_0}} J_{\Phi_M}^{(\epsilon, j)}(x) \ge I_j\bigl(x_0, \bar{\Gamma}_j\bigr), \quad j=0, 1, \ldots, n. \label{Eq20}
\end{align}
\end{lemma}

Let $\Gamma_j^{\circ}$ denote the interior of $\Gamma_j$ relative to $\partial \Omega$ and let $\bar{\Gamma}_j=\bar{\Gamma}_j^{\circ}$. Then, we have the following proposition.
\begin{proposition} \label{P1}
Suppose that Assumption~\ref{A1} holds, then, for $j=0, 1, \ldots, n$, we have
\begin{align}
 \epsilon\,\log \mathbb{P}_{x_0}^{(\epsilon, j)}\bigm(x^{\epsilon, j}(\tau_{j}) \in \Gamma_j\bigm) \rightarrow I_j\bigl(x_0, \Gamma_j\bigr) \quad \text{as} \quad \epsilon \rightarrow 0 \label{Eq21}
\end{align}
uniformly for all $x_0$ in any compact subset $\Lambda \subset \Omega$.
\end{proposition}

{\em Proof}: For any fixed $j \in \{0,1, \ldots, n\}$, it is suffices to show the following conditions
\begin{align}
  \limsup_{\epsilon \rightarrow 0} \epsilon\,\log \mathbb{P}_{x_0}^{(\epsilon, j)}\bigm(x^{\epsilon, j}(\tau_{j}) \in \Gamma_j\bigm) \le -I_j\bigl(x_0, \bar{\Gamma}_j\bigr), \label{Eq22}
\end{align}
and
\begin{align}
 \liminf_{\epsilon \rightarrow 0} \epsilon\,\log \mathbb{P}_{x_0}^{(\epsilon, j)}\bigm(x^{\epsilon, j}(\tau_{j}) \in \Gamma_j\bigm) \ge -I_j\bigl(x_0, \Gamma_j^{\circ}\bigr), \label{Eq23}
\end{align}
uniformly for $x_0 \in \Omega$ and for any $\Gamma_j \subset \partial \Omega$ with $\bar{\Gamma}_j=\bar{\Gamma}_j^{\circ}$.

Note that $I_j\bigl(x_0, \bar{\Gamma}_j^{\circ}\bigr) = I_j\bigl(x_0, \bar{\Gamma}_j\bigr)$ (cf. \eqref{Eq19}), then the upper bound in \eqref{Eq22} can be verified using the Ventcel-Freidlin estimate (see \cite[pp.\,332--334]{Fre06} or \cite{VenFre70}). 

On the other hand, to prove the lower bound in \eqref{Eq23}, we introduce a penalty function $\Phi_M$ (with $\Phi_M(y)=0$ for $y \in \Gamma$); and write $f^{(j)}(x)=f_M^{(j)}(x) \bigl(\equiv \mathbb{E}_{x_0}^{(\epsilon, j)} \bigl(\exp\bigl(-\frac{1} {\epsilon} \Phi_{M}\bigl(x\bigr) \bigr) \bigr)\bigr)$ and $J_{\Phi}^{(\epsilon, j)}=J_{\Phi_M}^{(\epsilon, j)}(x)$, with $\Phi_j=\Phi_M$. Then, from \eqref{Eq17}, we have
\begin{align}
 q^{(\epsilon,j)}\bigl(x_0, \Gamma_j\bigr) \le f_M^{(j)}(x_0), \label{Eq24}
\end{align}
for each $M$. Hence, using Lemma~\ref{L1} and noting further $J_{\Phi_M}^{(\epsilon, j)}(x_0) \ge I_j(x_0, \Gamma_j^{\circ})$, the lower bound in \eqref{Eq23} holds uniformly for all $x_0 \in \Lambda$. This completes the proof. \qed

In the following, using Proposition~\ref{P1}, we provide additional results on the exit positions of the state-trajectories $x^{\epsilon, j}(t)$ through the portion of the boundary $\Gamma_j$ for $j=0, 1, \dots, n$. 

For $x, y \in \bar{\Omega}$, we consider the following
\begin{align}
 I_j(\varphi, \tau_{j}) = \inf_{\varphi \in \mathcal{H}^{1}[0, T], \tau_{j} \ge 0} \int_{0}^{\tau_{j} \wedge T} L_j\big(\varphi(t),\dot{\varphi}(t)\big)dt,\quad \forall {j \in \mathcal{N}\cup \{0\}}, \label{Eq25}
\end{align}
where the infimum is taken among all $\varphi \in \mathcal{H}^{1}[0, T]$ and $\tau_{j} \ge 0$ such that $\varphi(0)=x_0$, $\varphi(\tau_{j})=y$ and $\varphi(t) \in \bar{\Omega}$ for all $t \in [0, \tau_{j} \wedge T]$. Then, using \eqref{Eq17}, we have
\begin{align}
 I_j\bigl(x_0, \Gamma_j\bigr) &= \inf_{y \in \Gamma} I_j\bigl(x_0, y\bigr), \notag \\
                                                    &= \min_{y \in \bar{\Gamma}} I_j\bigl(x_0, y\bigr),\quad j=0, 1, \dots, n, \label{Eq26}
\end{align} 
for $x_0 \in \Omega$ and $\Gamma_j \in \partial \Omega$.

Next, we will assume, in addition to Assumption~\ref{A1}, the followings (cf. \cite[pp\,359--360]{Fre06}).

\begin{assumption} \label{A2} ~
\begin{enumerate} [(a)]
\item $\Bigl\langle\bigl(A + \sum\nolimits_{i=1}^n B_i K_i \bigr)y,\, \gamma(y) \Bigr\rangle < 0$ and $\Bigl\langle\bigl(A + \sum\nolimits_{i \neq j} B_i K_i \bigr)y,\, \gamma(y) \Bigr\rangle < 0$ for $j=1, 2, \ldots, n$, where $\gamma(y)$ is the unit outward normal to $\Omega$ at $y \in \partial \Omega$.
\item For all $j \in \{0,1, \ldots, n\}$, let there exist a compact subset $\Lambda \subset \Omega$ such that:
\begin{enumerate} [(i)]
\item $I_j(x, y)=0$, $\forall x, y \in \Lambda$.
\item Let $\Lambda_{\delta}$ denote the $\delta$-neighborhood of $\Lambda$, and $\Omega_{\delta} = \Omega \setminus \bar{\Lambda}_{\delta}$. Then, there exists $c_{\delta}$ that tends to zero as $\delta \rightarrow 0$ such that
\begin{align}
 I_j^{\Omega_{\delta}}\bigl(x, y\bigr) \le I_j\bigl(x, y\bigr) + c_{\delta}, \quad \forall x, y \in \Omega \setminus \Lambda_{2\delta}, \label{Eq27}
\end{align}
where the minimum functional $I_j^{\Omega_{\delta}}$ is with respect to $\Omega_{\delta}$ (cf. equation~\eqref{Eq7}). 
\end{enumerate}
\end{enumerate}
\end{assumption}

Notice that the statements in Assumption~\ref{A2}(b) imply the following
\begin{align}
 I_j\bigl(x_1, y\bigr) = I_j\bigl(x_2, y\bigr), \quad \forall x_1, x_2 \in \Lambda, \quad j=0,1, \ldots, n. \label{Eq28}
\end{align}

Hence, for each $j = 0,1, \ldots, n$, if we let
\begin{align}
 V_j(x, \partial \Omega) = \inf_{y \in \partial \Omega} I_j\bigr(x, y\bigl), \quad x \in \Lambda, \label{Eq29}
\end{align}
and
\begin{align}
 \Sigma_j = \Bigl\{ y \in \partial \Omega \, \Bigl\vert  \, I_j\bigl(x, y\bigr) = V_j \bigl(x, \partial \Omega\bigr), \,\, x \in \Lambda \Bigr\}. \label{Eq30}
\end{align}
Then, we immediately obtain the following proposition.
\begin{proposition} \label{P2}
Suppose that Assumptions~\ref{A1} and \ref{A2} hold. Then, for any $\delta > 0$, $\operatorname{dist}\bigl(x^{\epsilon, j}, \Sigma_j \bigr) \rightarrow 0$ in probability as $\epsilon \rightarrow 0$ for each $j=0,1, \ldots, n$.
\end{proposition}

{\em Proof}: For any fixed $j \in \{0,1, \ldots, n\}$, let $S$ be open, with smooth boundary, and $\Lambda \subset S \subset \Lambda_{\delta}$ (where $\delta > 0$ is small enough such that $\bar{\Lambda}_{2\delta} \subset \Omega$). Further, let $\Omega_{\neg \bar{S}} = \Omega \setminus \bar{S}$ and let $\Gamma_j \subset \partial \Omega$ be closed with $ \Sigma_j \subset \bar{\Gamma}_j^{\circ}$ and $\bar{\Gamma}_j^{\circ} = \Gamma_j$. Then, for any $x_0 \in \Lambda$, there exits $\kappa > 0$ such that
\begin{align}
 I_j\bigl(x_0, \Gamma_j\bigr) =  V_j\bigl(x_0, \Sigma_j\bigr),  \label{Eq31}
\end{align}
and
\begin{align}
 I_j\bigl(x_0, \Gamma_j^{c}\bigr) =  V_j\bigl(x_0, \Sigma_j\bigr) + 2 \kappa, \label{Eq32}
\end{align}
where $\Gamma_j^{c} = \partial \Omega \setminus \Gamma_j$ for each $j=0,1, \ldots, n$.

Note that, from Assumption~\ref{A2}(b), one can choose small $\delta > 0$ such that
\begin{align}
 \max_{z \in \partial \Lambda_{2\delta}} I_j^{\Omega_{\neg \bar{S}}}\bigl(z, \Gamma_j\bigr) < V_j\bigl(z, \Sigma_j\bigr) + \kappa <  \min_{z \in \partial \Lambda_{2\delta}} I_j^{\Omega_{\neg \bar{S}}}\bigl(z, \Gamma_j^{c}\bigr), \label{Eq33}
 \end{align}
where the minimum functional $I_j^{\Omega_{\neg \bar{S}}}$ is with respect to $\Omega_{\neg \bar{S}}$. 

Then, from Proposition~\ref{P1}, we have the following
\begin{align}
 \lim_{\epsilon \rightarrow 0}\frac{q_{\neg \bar{S}}^{(\epsilon, j)}\bigl(z, \Gamma_j^{c}\bigr)}{q_{\neg \bar{S}}^{(\epsilon, j)}\bigl(z, \Gamma_j\bigr)} = 0, \quad j=0,1, \ldots, n, \label{Eq34}
\end{align}
uniformly for all $z \in \partial \Lambda_{2\delta}$. 

For $x_0 \in \Omega$ (with $x^{\epsilon, j}(0) = x_0$), let us define the following random time processes
\begin{align*}
 \tau_{j,0} &= \inf \Bigl \{ t \,\Bigl \vert \, x^{\epsilon, j}(t) \in \Omega_{\neg \bar{S}} \Bigr\},\\
 s_{j,k} &= \inf \Bigl \{ t \,\Bigl \vert \, t > \tau_{j, k-1}, \,\, k \ge 1, \,\, x^{\epsilon, j}(t) \in \Omega_{\neg \bar{S}} \Bigr\},\\
 \tau_{j,k} &= \inf \Bigl \{ t \,\Bigl \vert \, t > s_{j,k}, \,\, k \ge 1, \,\, x^{\epsilon, j}(t) \in \Omega_{\neg \bar{S}} \Bigr\}.
\end{align*}
Next, we consider the following events
\begin{align*}
 \mathscr{A}_{k}^{(j)} &= \Bigl \{ \tau_j = \tau_{j,k}, \,\, x^{\epsilon, j}(\tau_j) \in \Gamma_j \Bigr\},
\end{align*}
and 
\begin{align*}
 \mathscr{B}_{k}^{(j)} &= \Bigl \{ \tau_j = \tau_{j,k}, \,\, x^{\epsilon, j}(\tau_j) \in \Gamma_j^{c} \Bigr\}.
\end{align*}
Then, from the strong Markov property, we have
\begin{align}
 \mathbb{P}_{x_0}^{(\epsilon, j)}(\mathscr{A}_{k}^{(j)})= \mathbb{E}_{x_0}^{(\epsilon, j)} \Bigl (\chi_{\tau_j > s_{j,k}} q_{\neg \bar{S}}^{(\epsilon, j)}\bigl(x^{\epsilon, j}(s_{j,k}), \Gamma_j\bigr) \Bigr), \label{Eq35}
\end{align}
and 
\begin{align}
 \mathbb{P}_{x_0}^{(\epsilon, j)}(\mathscr{B}_{k}^{(j)})= \mathbb{E}_{x_0}^{(\epsilon, j)} \Bigl (\chi_{\tau_j > s_{j,k}} q_{\neg \bar{S}}^{(\epsilon, j)}\bigl(x^{\epsilon, j}(s_{j,k}), \Gamma_j^{c}\bigr) \Bigr), \label{Eq36}
\end{align}
where $\chi_{\tau_j > s_{j,k}}$ is an indicator function for the random event $\tau_j > s_{j,k}$ (with $k \ge 1$).\footnote{\label{FN6}Note that $\mathbb{P}_{x_0}^{(\epsilon, j)} \bigl(\mathscr{A}_{0}^{(j)} \bigcup \mathscr{B}_{0}^{(j)}\bigl) \rightarrow 0$ in probability as $\epsilon \rightarrow 0$ for each $j=0, 1, \ldots, n$.} Note that, from \eqref{Eq35}, for any $\ell >0$, there exits an $\epsilon_{\ell}>0$ such that  
\begin{align}
 q_{\neg \bar{S}}^{(\epsilon, j)}(x_0, \Gamma_j^{c}) \le \ell\, q_{\neg \bar{S}}^{(\epsilon, j)}(x_0, \Gamma), \quad \forall z \in \partial \Lambda_{2\delta},\,\, \forall  \epsilon \in (0, \epsilon_{\ell}). \label{Eq37}
\end{align}
Since $x^{\epsilon, j}(s_{j,k}) \in \partial \Lambda_{2\delta}$, then we have $\mathbb{P}_{x_0}^{(\epsilon, j)} \Bigl(\mathscr{B}_{k}^{(j)}\Bigr) \le \ell \, \mathbb{P}_{x_0}^{(\epsilon, j)} \Bigl(\mathscr{A}_{k}^{(j)}\Bigr)$. Moreover, we have
\begin{align}
\sum\nolimits_{k} \mathbb{P}_{x_0}^{(\epsilon, j)} \Bigl(\mathscr{A}_{k}^{(j)}\Bigr) &\le \sum\nolimits_{k} \mathbb{P}_{x_0}^{(\epsilon, j)} \Bigl(\mathscr{A}_{k}^{(j)} \bigcup \mathscr{B}_{k}^{(j)}\Bigr), \notag\\
                                                                            &= \mathbb{P}_{x_0}^{(\epsilon, j)} \Bigl(\tau_j < \infty \Bigr) \biggl(\equiv 1, \,\, \text{for each} \,\, j=0,1, \ldots, n\biggr). \label{Eq38}                                                                         
\end{align}
Hence, for $\epsilon \in (0, \epsilon_{\ell})$, we have (see also Footnote~\ref{FN6})
\begin{align}
 \mathbb{P}_{x_0}^{(\epsilon, j)} \Bigl(x^{\epsilon, j}(\tau_j) \in \Gamma_j^{c} \Bigr) &= \sum\nolimits_{k} \mathbb{P}_{x_0}^{(\epsilon, j)} \Bigl(\mathscr{B}_{k}^{(j)}\Bigr), \notag \\
                                                 & \le \mathbb{P}_{x_0}^{(\epsilon, j)} \Bigl(\mathscr{B}_{0}^{(j)}\Bigr) + \ell,  \quad j=0,1, \ldots, n. \label{Eq39}
\end{align}
Since $\ell$ is arbitrary, this completes the proof. 
\qed

Note that the above proposition, i.e., Proposition~\ref{P2}, is connected to the boundary value problem, when one is also interested on the position of state-trajectories at the first time of their exit from the boundary $\Gamma_j \subset \partial\Omega$. For the boundary value problem of Section~\ref{S2} (cf. \cite[pp.\,371--272]{Fre06}), with $\mathcal{L}^{(\epsilon, j)} f_M^{(j)}(x) = 0$ in $\Omega$ and $f_M^{(j)}(x) = \mathbb{E}_{x_0}^{(\epsilon, j)} \bigl(\exp\bigl(-\frac{1} {\epsilon} \Phi_M\bigl(x\bigr) \bigr) \bigr)$ on $\partial\Omega$. For example, if $\Sigma_j$ consists of a single point, say $y^{\ast}$, then $f_M^{(j)}(x) \rightarrow \mathbb{E}_{x_0}^{(\epsilon, j)} \bigl(\exp\bigl(-\frac{1} {\epsilon} \Phi_M\bigl(y^{\ast}\bigr) \bigr) \bigr)$ as $\epsilon \rightarrow 0$ for all $x_0 \in \Omega$. Moreover, if $\Lambda$ consists of a single point $x^{\ast} \in \Omega$, then we have the following
\begin{align}
 - \lim_{\epsilon \rightarrow 0} \epsilon \log \mathbb{E}_{x^{\ast}}^{(\epsilon, j)} \Bigl(x^{\epsilon, j}(\tau_{j}) \in \Gamma_j\Bigr) &= \min_{y \in \partial \Omega} I_j(x^{\ast}, y), \notag \\
                                                                                                      &= V_j\bigl(x^{\ast}, \Sigma_j\bigr), \quad j=0,1, \ldots, n, \label{Eq40}
\end{align}
which is equivalent to the result of Proposition~\ref{P1}.\footnote{Note that the last portion of the state-trajectories, prior reaching the boundary $\Gamma_j$ lies in the neighborhood of $\varphi(t) \in \mathcal{H}^{1}[0, T]$ for which $I_j(x^{\ast}, \Gamma_j)$ differs little from $V_j\bigl(x^{\ast}, \Sigma_j\bigr)$ for each $ j=0,1, \ldots, n$.}

\section{Further remarks} \label{S5}
In this section, we briefly comment on the implication of our results on a co-design problem -- when one is also interested in either evaluating the performance or finding a set of sub-optimal redundant controllers for the multi-channel system, while estimating the exit probabilities or the asymptotic bounds on the mean exit-time of the state-trajectories from the domain $\Omega$.\footnote{Note that, in Section~\ref{S4}, we provide estimates for the exit probabilities of state-trajectories from the boundary of the given domain $\Omega$ for each particular operating condition (i.e., during the nominal operating condition or any single-channel failure in the multi-channel system).} 

In particular, here we outline a multi-objective embedded optimization framework -- where the problem of optimal exit probabilities and finding a set of stabilizing feedbacks for the multi-channel system can be considered as a composite goal-oriented optimization problem (see, e.g., \cite{GemHa75}). Hence, the composite optimization problem (which also embeds additional  subproblems) can be reformulated as follows 
\begin{align}
\left.\begin{array}{l}
 \min\boldsymbol{\gamma}  \\ 
  \text{subject to} \\
 \quad I_i\bigl(x_0, \Gamma_i\bigr) - \gamma_i w_i \le I_0\bigr(\varphi^{\ast}, {\tau_0^{\ast}}, K^{\ast}\bigr), \,\, \text{with} \,\, \Gamma_i \subset \partial \Omega \,\, \& \,\,  K^{\ast} \in \mathcal{K} \\
 \quad  x_0 \in \Omega  \quad \text{(initial condition)} \\
 \quad w_i > 0, \quad \text{with} \quad \sum\nolimits_{i =1}^n w_i\,= 1  \\
 \quad \gamma_i \quad \text{(unrestricted scalar variables)} 
\end{array} \right\} \label{Eq41}
\end{align}
where $w_i$'s are the weighting factors and the vector $\boldsymbol{\gamma}$ is given by $[\gamma_1, \gamma_2, \ldots, \gamma_n]^T$. Moreover, $I_0(\varphi^{\ast}, \tau_0^{\ast}, K^{\ast})$, which corresponds to the nominal operating condition (i.e., without any fault in the system), is given by  
\begin{align}
 I_0\bigl(\varphi^{\ast}, {\tau_0^{\ast}}, K^{\ast}\bigr) = \sup_{K \in \mathcal{K}} \, \inf_{\varphi \in \mathcal{H}^{1}[0, T], \tau_0 \ge 0} \int_{0}^{\tau_0} L_0\big(\varphi(t),\dot{\varphi}(t)\big)dt, \label{Eq42}
\end{align}
where the $I_i\bigl(x_0, \Gamma_i\bigr)$'s (together with the boundary conditions $x^{\epsilon, i}(\tau_{i}) \in \Gamma_i$ for $i=1,2, \ldots, n$) are assumed to satisfy the optimization subproblems in \eqref{Eq26} (cf. equation~\eqref{Eq40}). Note that such class of stabilizing state-feedbacks can further be restricted to satisfy additional assumptions. Here, we remark that the {\em max-min problem} in \eqref{Eq42} for the exit probabilities has been studied in the past (see, e.g., \cite{FleTs81} in the context of differential games for a general admissible class of controls; and see also \cite{EvaIsh85} or \cite{Fle78} via viscosity solution techniques).

\begin{remark} \label{R4}
Finally, we remark that the composite optimization problem in \eqref{Eq41} is useful for selecting the most appropriate $n$-tuple of stabilizing state-feedbacks from the set $\bigl\{K^{\ast}\bigr\}_{\nu} \in \mathcal{K}$ that confines the state-trajectories $x^{\epsilon,j}(t)$ to the prescribed domain $D$ for a certain duration, while the system performances are, in some sense, associated with those portions of the boundary of the given domain (i.e., $\Gamma_j \subset \partial \Omega$, $j=0,1, \ldots, n$). Such a composite optimization problem, although computationally demanding, can be sub-optimally solved by relaxing some of the constraints.
\end{remark}

%
%
%
%
%

\end{document}